\documentclass[preprint,3p,10pt,authoryear]{elsarticle}

\usepackage{amssymb,amsthm}
\usepackage{lineno}  

\usepackage{graphicx,mathtools,natbib,caption,geometry}
\usepackage{hyperref}
\hypersetup{
     colorlinks   = true,
     citecolor    = blue
}

\usepackage{listings}
\usepackage{appendix}

\newtheorem{thm}{Theorem}

\newtheorem{rem}{Remark}
\newtheorem{examp}[rem]{Example}
\newtheorem{condition}{Condition}

\newenvironment{proofff}[1][Proof]{\begin{trivlist}
\item[\hskip \labelsep {\bfseries #1}]}{\end{trivlist}}

\makeatletter \@addtoreset{equation}{section} \makeatother

\newcommand{\N}{\mathbb{N}}
\newcommand{\R}{\mathbb{R}}
\newcommand{\PP}{\mathbb{P}}
\newcommand{\E}{\mathbb{E}}
\newcommand{\Var}{\mathbb{V}\mathrm{ar}}
\newcommand{\leqdef}{\vcentcolon=}
\newcommand{\rd}{{\rm d}}
\newcommand{\ind}{\mathbb{I}}
\newcommand{\ii}{\mathrm{i}}
\newcommand{\bb}[1]{\boldsymbol{#1}}

\allowdisplaybreaks

\makeatletter
\def\ps@pprintTitle{%
 \let\@oddhead\@empty
 \let\@evenhead\@empty
 \def\@oddfoot{\centerline{\thepage}}%
 \let\@evenfoot\@oddfoot}
\makeatother

\begin{document}

\begin{frontmatter}

    \title{A counterexample to the central limit theorem for pairwise independent random variables having a common arbitrary margin}%

    \author[a1]{Benjamin Avanzi\fnref{fn1}} \ead{b.avanzi@unimelb.edu.au}
    \author[a2]{Guillaume Boglioni Beaulieu\corref{cor2}\fnref{fn2}} \ead{g.boglionibeaulieu@unsw.edu.au}
    \author[a3]{Pierre Lafaye de Micheaux} \ead{lafaye@unsw.edu.au}
    \author[a4]{Fr\'ed\'eric Ouimet\fnref{fn4}} \ead{ouimetfr@caltech.edu}
    \author[a2]{Bernard Wong\fnref{fn1}} \ead{bernard.wong@unsw.edu.au}

    \address[a1]{Centre for Actuarial Studies, Department of Economics, University of Melbourne, VIC 3010, Australia.}%
    \address[a2]{School of Risk and Actuarial Studies, UNSW Sydney, NSW 2052, Australia.}%
    \address[a3]{School of Mathematics and Statistics, UNSW Sydney, NSW 2052, Australia.}%
    \address[a4]{PMA Department, California Institute of Technology, CA 91125, Pasadena, USA.}%

    \cortext[cor2]{Corresponding author}%

    \fntext[fn1]{B.\ A.\ and B.\ W.\ are supported by Australian Research Council's Linkage (LP130100723) and Discovery (DP200101859) Projects funding schemes.}%
    \fntext[fn2]{G.\ B.\ B.\ acknowledges financial support from UNSW Sydney under a University International Postgraduate Award, from UNSW Business School under a supplementary scholarship, and from the FRQNT (B2).}%
    \fntext[fn4]{F.\ O.\ is supported by a postdoctoral fellowship from the NSERC (PDF) and the FRQNT (B3X supplement).}%

    \begin{abstract}
        The Central Limit Theorem (CLT) is one of the most fundamental results in statistics. It states that the standardized sample mean of a sequence of $n$ \emph{mutually} independent and identically distributed random variables with finite first and second moments converges in distribution to a standard Gaussian as $n$ goes to infinity. In particular, \emph{pairwise} independence of the sequence is generally not sufficient for the theorem to hold. We construct explicitly a sequence of pairwise independent random variables having a common but arbitrary marginal distribution $F$ (satisfying very mild conditions) for which the CLT is not verified. We study the extent of this `failure' of the CLT by obtaining, in closed form, the asymptotic distribution of the sample mean  of our sequence. This is illustrated through several theoretical examples, for which we provide associated computing codes in the \texttt{R} language.
    \end{abstract}


\end{frontmatter}

\section{Introduction}

The aim of this paper is to construct explicitly a sequence of {\it pairwise} independent and identically distributed (p.i.i.d.) random variables (r.v.s)\ whose common margin $F$ can be chosen arbitrarily (under very mild conditions) and for which the (standardized) sample mean is {\it not} asymptotically Gaussian. We give a closed-form expression for the limiting distribution of this sample mean. It is, to the best of our knowledge, the first example of this kind for which the asymptotic distribution of the sample mean is explicitly given, and known to be skewed and heavier tailed than a Gaussian distribution, for {\it any} choice of margin. Our sequence thus illustrates nicely why {\it mutual} independence is such a crucial assumption in the Central Limit Theorem (CLT). Furthermore, it allows us to quantify how far away from the Gaussian distribution one can get under the less restrictive assumption of {\it pairwise} independence.

Recall that the classical CLT is stated for a sequence of i.i.d.\ random variables where the first `i' in the acronym stands for `independent', which itself stands for `{\it mutually} independent' (while the last `i.d.' stands for {\it identically distributed}). Now, it is known that \textit{pairwise} independence among random variables is a necessary \emph{but not sufficient} condition for them to be \textit{mutually} independent. The earliest counterexample can be attributed to \citet{bernstein1927}, followed by a few other authors, e.g., \citet{geisser1962, pierce1969, joffe1974, bretagnolle1995, derriennic2000}.
However, from these illustrative examples alone it can be hard to understand how bad of a substitute to mutual independence pairwise independence is. One way to study this question is to consider those fundamental theorems of mathematical statistics that rely on the former assumption; do they `fail' under the weaker assumption of \textit{pairwise} independence? A definite answer to that question is beyond the scope of this work, as it depends on which theorem is considered. Nevertheless, note that the Law of Large Numbers, even if almost always stated for mutually independent r.v.s, does hold under pairwise independence \citep{etemadi1981}. The same goes for the second Borel-Cantelli lemma, usually stated for mutually independent events but valid for pairwise independent events as well \citep{erdos1959}. The CLT, however, does `fail' under pairwise independence. Since it is arguably the most crucial result in all of statistics, we will concentrate on this theorem from now on.

The classical CLT is usually stated as follows. Given a sequence $\{X_j, j \geq 1\}$ of i.i.d.\ r.v.s with mean $\mu$ and standard deviation $0 < \sigma<\infty$, we have, as $n$ tends to infinity,
\begin{equation}
S_n \leqdef \frac{\sum_{j=1}^n X_j - \mu n}{\sigma\sqrt{n}}\stackrel{d}{\longrightarrow}Z,
\end{equation}
where the random variable $Z$ has a standard Gaussian distribution, noted thereafter $N(0,1)$, and `$\stackrel{d}{\longrightarrow}$' denotes convergence in distribution. \citet{revesz1965} were the first to provide a pairwise independent sequence for which $S_n$ does not converge in distribution to a $N(0,1)$. For their sequence, which is binary (i.e., two-state), $S_n$ converges to a standardized $\chi^2_1$ distribution. \citet[Example~5.45]{romano1986} provide a two-state, and \cite*{bradley1989} a three-state, pairwise independent sequence for which $S_n$ converges in probability to $0$. \cite*{janson1988} provides a broader class of pairwise independent `counterexamples' to the CLT, most defined with $X_j$'s having a continuous margin and for which $S_n$ converges in probability to $0$. The author also constructs a pairwise independent sequence of $N(0,1)$ r.v.s for which $S_n$ converges to the random variable $S = R \cdot Z$, with $R$ a r.v.\ whose distribution can be arbitrarily chosen among those with support $[0,1]$, and $Z$ a $N(0,1)$ r.v.\ independent of $R$. The r.v.\ $S$ can be seen as `better behaved' than a $N(0,1)$, in the sense that it is symmetric with a variance smaller than $1$ (regardless of the choice of $R$).
\citet*[][Section 2.3]{cuesta1991} construct a sequence $\{X_j, j \geq 1\}$ of r.v.s taking values uniformly on the integers $\{0,1,\ldots, p-1\}$, with $p$ a prime number, for which $S_n$ is `worse behaved' than a $N(0,1)$. Indeed, their $S_n$ converges in distribution to a mixture (with weights $(p-1)/p$ and $1/p$ respectively) of the constant $0$ and of a centered Gaussian r.v.\ with variance $p$. This distribution is symmetric but it has heavier tails than that of a $N(0,1)$.

Other authors go beyond \textit{pairwise} independence and study the CLT under `$K$-tuplewise independence', for $K \geq 3$. A random sequence $\{X_j, j \geq 1\}$ is said $K$-tuplewise independent if for every choice of $K$ distinct integers $j_1, \ldots, j_K$, the random variables $X_{j_1}, \ldots, X_{j_K}$ are mutually independent. \cite*{kantorovitz2007} provides an example of a triple-wise independent two-state sequence for which $S_n$ converges to a `misbehaved' distribution ---that of $Z_1 \cdot Z_2$, where $Z_1$ and $Z_2$ are independent $N(0,1)$. 
\cite*{pruss1998} presents a sequence of $K$-tuplewise independent random variables $\{X_j, j \geq 1 \}$ taking values in $\{-1,1\}$ for which the asymptotic distribution of $S_n$ is never Gaussian, whichever choice of $K$. \cite*{bradley2009} extend this construction to a \textit{strictly stationary} sequence of $K$-tuplewise independent r.v.s whose margin is uniform on the interval $[-\sqrt{3}, \sqrt{3}]$.
\cite*{weakley2013} further extends this construction by allowing the $X_j$'s to have any symmetrical distribution (with finite variance). 

In the body of research discussed above, a non-degenerate and explicit limiting distribution for $S_n$ is obtained only for very specific choices of margin for the $X_j$'s.
In this paper, we allow this margin to be almost any non-degenerate distribution, yet we still obtain explicitly the limiting distribution of $S_n$. This distribution depends on the choice of the margin, but it is always skewed and heavier tailed than that of a Gaussian. 

The rest of the paper is organized as follows. In Section~\ref{clt_example_setting}, we construct our pairwise independent sequence $\{X_j, j \geq 1 \}$. In Section~\ref{asymptotic-distribution-of-the-mean}, we derive explicitly the asymptotic distribution of the standardized average of that sequence, while in Section~\ref{properties-of-s}, we study key properties of such a distribution. In Section~\ref{conclusion_section}, we conclude.




\section{Construction of the Gaussian-Chi-squared pairwise independent sequence}\label{clt_example_setting}

In this section, we build a sequence $\{X_j, j \geq 1\}$ of p.i.i.d.\ r.v.s for which the CLT does \textit{not} hold. We show in Section~\ref{properties-of-s} that the asymptotic distribution of the (standardized) sample mean of this sequence can be conveniently written as that of the sum of a Gaussian r.v.\ and of an independent scaled   Chi-squared r.v. Importantly,  the r.v.s  forming this sequence have a common (but arbitrary) marginal distribution $F$ satisfying the following condition:
\begin{condition}\label{cond:F}
     For any r.v.\ $W\sim F$, the variance $\Var(W)$ is finite and there exists a Borel set $A$ for which $\PP(W\in A) = \ell^{-1}$, for some integer $\ell \geq 2$, and $\mathbb{E}[W | A] \neq \mathbb{E}[W | A^c]$.
\end{condition}

\noindent
As long as the variance is finite, the restriction on $F$ includes {\it all} distributions with an absolutely continuous part on some interval. It also includes {\it almost all} discrete distributions with at least one weight of the form $\ell^{-1}$; see Remark~\ref{rem:discretedist}. Also, note that, for a given $F$, many choices for $A$ (with possibly different values of $\ell$) could be available, depending on $F$.

We begin our construction of the sequence $\{X_j, j \geq 1\}$ by letting $F$ be a distribution satisfying Condition \ref{cond:F}. For a r.v.\ $W \sim F$,
let $A$ be any Borel set such that
\begin{equation}\label{eq:w.A}
\PP(W\in A) = \ell^{-1}, \quad \text{for some integer } \ell\geq 2.
\end{equation}
Then, for an integer $m \geq 2$, let $\bb{M}_1, \ldots, \bb{M}_m$ be a sequence of i.i.d.\  r.v.s defined on a common probability space $(\Omega, \mathcal{F}, \mathbb{P})$ having the following distribution:
\begin{equation}\label{eq:multinomial}
     \textup{Multinomial}\left( 1\, ; (p_1 = \ell^{-1},p_2=\ell^{-1},\dots,p_{\ell}=\ell^{-1})\right).
\end{equation}
For all pairs \((\bb{M}_i, \bb{M}_j)\), $1 \leq i < j \leq m$, define a r.v.\ \(D_{i,j}\) as
\begin{equation}
    D_{i,j} =
    \left\{\hspace{-1mm}
    \begin{array}{ll}
    1, & \textup{if } \bb{M}_i = \bb{M}_j, \\
    0, & \textup{otherwise}.
    \end{array}
    \right.
\end{equation}
The $D_{i,j}$ are p.i.i.d (but not mutually independent) with  $\mathbb{P}(D_{i,j} = 1)= \ell^{-1}$; see Remark~\ref{rem:threeconds}. For convenience, we refer to these $n = \binom{m}{2}$ random variables $D_{1,2}, D_{1,3}, \ldots, D_{1,m}, D_{2,3}, D_{2,4}, \ldots, D_{m-1,m}$ simply as
\begin{align}\label{the_D_sequence}
    D_1, \ldots, D_n,
\end{align}
where for $1 \leq i < j \leq m$, $D_{k(i,j)} \leqdef D_{i,j}$ with $k(i,j) = [i(2m-1)-i^2]/2+j-m.$
Note that when $\ell = 2$ and $p_1 = p_2 = 1/2$, the $\bb{M}_j$'s can be identified simply as $\text{Bernoulli}(1/2)$ r.v.s, and the sequence \eqref{the_D_sequence} is equivalent to a pairwise independent sequence first mentioned in \citet{geisser1962} and for which we already know that the CLT does not hold  \citep[see][]{revesz1965}.


From the sequence $D_1, \ldots, D_n$, we now  construct a new pairwise independent sequence $X_1, \ldots, X_n$ such that $X_k \sim F$ for all $k = 1, \ldots, n$. Define $U$ and $V$ to be the truncated versions of $W$, respectively off and on the set $A$:
\begin{equation}
    U \stackrel{d}{=} W \mid \{W\in A^c\}, \qquad V \stackrel{d}{=} W \mid \{W\in A\},
\end{equation}
and denote
\begin{align}\label{def:mu_U}
    \mu_U \coloneqq \E[U], \qquad \mu_V \coloneqq \E[V].
\end{align}
Then, consider $n$ independent copies of $U$, and independently $n$ independent copies of $V$:
\begin{equation}
    U_1, \ldots, U_n, \stackrel{\textup{i.i.d.}}{\sim} F_U, \qquad V_1, \ldots, V_n \stackrel{\textup{i.i.d.}}{\sim} F_V,
\end{equation}
both defined on the probability space $(\Omega, \mathcal{F}, \mathbb{P})$.
Finally, for $\omega \in \Omega$ and for  $k = 1,\ldots, n$, construct
\begin{equation}\label{the_X_sequence}
    X_k(\omega) = \left\{\hspace{-1mm}
    \begin{array}{ll}
    U_k (\omega), &  \textup{if } D_k(\omega) = 0, \\
    V_k (\omega), &  \textup{if } D_k(\omega) = 1.
    \end{array}
    \right.
\end{equation}
By conditioning on $D_k$, one can check  that
\begin{align}\label{eq:F_Xk}
    F_{X_k}(x)
    &= (1 - \ell^{-1}) F_{U_k}(x) + \ell^{-1} F_{V_k}(x) = F(x).
\end{align}
In the next section, we will derive the asymptotic distribution of the sample mean of those $X$'s, and see that it is \textit{not} Gaussian.

The failure of the CLT for this sequence \eqref{the_D_sequence} can be explained heuristically as follows. Within the sequence $D_1, \ldots, D_n$, there can be a `very high' proportion of $1$'s. This occurs if the sequence $\bb{M}_1, \ldots, \bb{M}_n$ contains a large proportion of equal vectors. However, by definition of $D_{i,j}$, to have a very large proportion of $0$'s among the $D$'s, one would require a large proportion of pairs $(\bb{M}_i,\bb{M}_j), 1 \leq i < j\leq m$ to be such that $\bb{M}_i \neq \bb{M}_j$. This is impossible, since \textit{all} the possible pairs are used to form the sequence of $D$'s. This very asymmetrical situation makes the asymptotic distribution of the standardized sample mean of the $D$'s 
highly skewed to the right. 

\begin{rem}\label{rem:threeconds}
In Condition \ref{cond:F}, the restriction $\PP(W\in A) = \ell^{-1}$ for some integer $\ell$ may seem  arbitrary. Likewise, in \eqref{eq:multinomial} the choice $p_i = \ell^{-1}$ for $i=1, \ldots, \ell$ may also seem arbitrary. We establish here that none of these choices are arbitrary. Indeed, assume first that the only restriction on $p_1, p_2, \ldots, p_{\ell}\in (0,1)$ is that
\begin{equation}\label{eq:p1.p2.condition}
    \begin{aligned}
    &(1) ~:~ p_1 + p_2 + \dots + p_{\ell} = 1, \\
    &(2) ~:~ p_1^2 + p_2^2 + \dots + p_{\ell}^2 = w, \\
    &(3) ~:~ p_1^3 + p_2^3 + \dots + p_{\ell}^3 = w^2,
    \end{aligned}
\end{equation}
for some $w \in (0,1)$. Condition $(1)$ is necessary for the multinomial in \eqref{eq:multinomial} to be well-defined, and conditions $(2)$ and $(3)$ are rewritings of
\begin{equation}\label{eq:D.condition.1}
    \PP(D_{i,j} = 1) = w, ~~1 \leq i < j \leq m,
\end{equation}
and
\begin{equation}\label{eq:D.condition.2}
    \PP(D_{i,j} = 1, D_{j,k} = 1) = \PP(D_{i,j} = 1) \PP(D_{j,k} = 1), \quad 1 \leq i < j < k \leq m,
\end{equation}
which are sufficient to guarantee that the $D$'s are respectively identically distributed and pairwise independent. Now, the solution $p_i = \ell^{-1}$ to \eqref{eq:p1.p2.condition} is \textit{unique}. Indeed, by squaring  condition~(2) in~\eqref{eq:p1.p2.condition} then applying  the Cauchy-Schwarz inequality, one gets
\begin{equation}\label{eq:system.solve}
    w^2 = \Big(\sum_{i=1}^{\ell} p_i^{3/2} p_i^{1/2}\Big)^2 \leq \sum_{i=1}^{\ell} p_i^3 \sum_{i=1}^{\ell} p_i = \sum_{i=1}^{\ell} p_i^3
\end{equation}
where the last equality comes from condition~(1) in~\eqref{eq:p1.p2.condition}.
Then, condition $(3)$ requires that we have the equality in \eqref{eq:system.solve}, and this happens if and only if $p_i^{3/2} = \lambda p_i^{1/2}$ for all $i\in\{1,\ldots,n\}$ and for  some $\lambda\in \R$. In turn, this implies $p_i = \lambda = \ell^{-1}$ because of $(1)$ and since $p_i > 0$, which then implies $w  = \ell^{-1}$ by $(2)$.
This reasoning shows that we cannot extend our method to an arbitrary $\PP(W\in A)\in (0,1)$ in \eqref{eq:w.A}. 
\end{rem}

\begin{rem}\label{rem:discretedist}
There is no easy characterization of all discrete distributions with finite variance such that $\PP(W\in A) = \ell^{-1}$ for some Borel set $A$ and some integer $\ell \geq 2$, but for which the last part of Condition \ref{cond:F} is not satisfied. However, the proportion of such distributions can be expected to be very small. As a simple but convincing example, consider the set of discrete distributions on three values $-\infty < x < y < z < \infty$ with weights $p_x, p_y, p_z\in (0,1)$. The variance is finite and say one of the three $p$'s has the form $\ell^{-1}$ for some integer $\ell^{-1}$. The only way that $\mathbb{E}[W | A] = \mathbb{E}[W | A^c]$ is satisfied is by having $A$ contain $y$ and only $y$ so that we must have $p_y = \ell^{-1}$, $p_x = p$ and $p_z = (1 - p - \ell^{-1})$ for some parameter $p\in (0,1)$, and $x p_x + z p_z = y \, (1 - \ell^{-1})$. In other words, once $\ell$ is fixed, there is only freedom in the choice of $x$, $z$ and $p$. If we remove the restriction $\mathbb{E}[W | A] = \mathbb{E}[W | A^c]$ (i.e.\ $x p_x + z p_z = y \, (1 - \ell^{-1})$), it gives us at least one more dimension of freedom in the selection of $x, y, z, p_x, p_y, p_z$. Hence, in this case, the proportion is actually ``$0$''. An analogous argument can be made for other discrete distributions of this kind. The restriction $\mathbb{E}[W | A] = \mathbb{E}[W | A^c]$ will always remove a dimension of freedom in the choice of the range of values or the weights.
\end{rem}



\section{Main result}\label{asymptotic-distribution-of-the-mean}

We now state our main result.

\begin{thm}\label{limiting_Sn_thm}
    Let $X_1, \ldots, X_n$ be random variables
    defined as in \eqref{the_X_sequence} and denote their mean and variance by $\mu$ and $\sigma^2$, respectively. Then under Condition~\ref{cond:F},

    \begin{itemize}
        \item[(a)] \textit{$X_1, \ldots, X_n$ are pairwise independent;}
        \item[(b)] As $m\to \infty$ (and hence as $n\to \infty$), the standardized sample mean $S_n \leqdef \big(\sum_{k=1}^n X_k - \mu n\big) / \sigma \sqrt{n}$ converges in distribution to a random variable
            \begin{equation}\label{eq:limit.S}
                S \leqdef \sqrt{1 - r^2} Z + r \, \chi,
            \end{equation}
            where $Z\sim N(0,1)$, $\chi$ is independently  distributed as a standardized $\chi_{\ell-1}^2$  and $r \leqdef \sqrt{\ell^{-1} (1 - \ell^{-1})} (\mu_V - \mu_U)/\sigma$ with $\mu_U, \mu_V$ defined in \eqref{def:mu_U}.
    \end{itemize}
\end{thm}

\begin{rem}
Interestingly, since a standardized chi-squared distribution converges to a standard Gaussian as its degree of freedom tends to infinity, we see that $S\stackrel{d}{\longrightarrow} N(0,1)$ as $\ell\to \infty$.
\end{rem}

\begin{rem}
When removing the restriction $\mathbb{E}[W | A] \neq \mathbb{E}[W | A^c]$ in Condition \ref{cond:F}, the case $r=0$ (i.e.\ $\mu_U = \mu_V$) is possible, so our construction also provides a new instance of a pairwise independent (but not mutually independent) sequence for which the CLT does hold.
\end{rem}

\begin{proofff}[Proof of Theorem \ref{limiting_Sn_thm}.]
    Proving (a) is straightforward. Simple calculations show that \(D_1, \ldots, D_n\) are pairwise independent; recall \eqref{eq:D.condition.2}. Now, for any $k, k' \in \{1,2,\ldots,n\}$ with $k \neq k'$, the r.v.s
    \begin{equation*}
        D_k, U_k, V_k, D_{k'}, U_{k'}, V_{k'}
    \end{equation*}
    are mutually independent and one can write $X_k = g(D_k, U_k, V_k)$ and $X_{k'} = g(D_{k'}, U_{k'}, V_{k'})$, for $g$ a Borel-measurable function. Since $X_k$ and $X_{k'}$ are integrable, the result follows from \citet[][Section 4.1, Corollary 2]{pollard2002}.

    The proof of $(b)$ is more involved.
    We prove \eqref{eq:limit.S} by obtaining the limit of the characteristic function of $S_n$, and then by invoking L{\'e}vy's continuity theorem.
    Namely, we show that, for all $t\in \R$,
    \begin{align}\label{cf_S}
        \varphi_{S_n}(t) \underset{m\to \infty}{\longrightarrow} \varphi_{\sqrt{1 - r^2} Z}(t) \cdot \varphi_{r \, \chi}(t) = e^{-\frac{1}{2} (1 - r^2) t^2} \hspace{-0.5mm}\cdot e^{-\ii t r \sqrt{(\ell - 1) / 2}} \Big(1 - \ii t r \sqrt{2 / (\ell - 1)}\Big)^{-(\ell-1)/2}.
    \end{align}
    First, let us define by $N_i = N_i(m)$ the number of $\bb{M}_j$'s in the $i$-th category, $i = 1,2,\dots,\ell$, within the multinomial sample $\{\bb{M}_j;~j=1,\ldots,m\}$.
    Then, $\bb{N} \leqdef (N_1,\dots,N_{\ell}) \sim \textup{Multinomial}\hspace{0.3mm}(m, (\ell^{-1}, \dots, \ell^{-1}))$.
    Importantly, if $\bb{N}$ is known, then the number $p(\bb{N})$ of $1$'s in the sequence $\{D_j;~ j=1, \ldots,n\}$ can be deduced as
    \begin{align}\label{np}
        p(\bb{N})
        &= \sum_{i=1}^{\ell - 1} \binom{N_i}{2} \ind_{\{N_i \geq 2\}} + \binom{m - \sum_{i=1}^{\ell-1} N_i}{2} \ind_{\{m - \sum_{i=1}^{\ell-1} N_i \geq 2\}} \notag \\[1mm]
        &= \sum_{i=1}^{\ell-1} \frac{N_i (N_i - 1)}{2} + \frac{(m - \sum_{i=1}^{\ell-1} N_i) (m - \sum_{i=1}^{\ell-1} N_i - 1)}{2} \notag \\[1mm]
        &= \frac{1}{2} \sum_{i=1}^{\ell-1} N_i^2 + \frac{1}{2} \sum_{i=1}^{\ell-1} \sum_{i'=1}^{\ell-1} N_i N_{i'} - m \sum_{i=1}^{\ell-1} N_i + \frac{m(m-1)}{2} \notag \\[1mm]
        &= \frac{1}{2} \sum_{i=1}^{\ell-1} (N_i - m \ell^{-1})^2 + \frac{1}{2} \sum_{i=1}^{\ell-1} \sum_{i'=1}^{\ell-1} (N_i - m \ell^{-1}) (N_{i'} - m \ell^{-1}) - \frac{\ell (\ell - 1)}{2} m^2 \ell^{-2} + \frac{m(m-1)}{2} \notag \\[1mm]
        &= \frac{m \ell^{-1}}{2} \sum_{i=1}^{\ell-1} \sum_{i'=1}^{\ell-1} \Big(\frac{1}{\ell^{-1}} \ind_{\{i = i'\}} + \frac{1}{\ell^{-1}}\Big) \frac{\big(N_i - m \ell^{-1}\big)}{\sqrt{m}} \frac{\big(N_{i'} - m \ell^{-1}\big)}{\sqrt{m}} - \frac{m}{2} (1 - \ell^{-1}) + n \ell^{-1},
    \end{align}
    where $\ind_B$ denotes the indicator function on the set $B$.
    The covariances of a $\textup{Multinomial}\hspace{0.3mm}(m, (p_1,p_2,\dots,p_{\ell}))$ distribution are well known to be $m \Sigma$ where $\Sigma_{i,i'} = p_i \ind_{\{i=i'\}} - p_i p_{i'}$, for $1 \leq i,i' \leq \ell-1$, and it is also known that $(\Sigma^{-1})_{i,i'} = p_i^{-1} \ind_{\{i=i'\}} + p_{\ell}^{-1}, ~ 1 \leq i,i' \leq \ell-1$; see \cite[eq.\hspace{0.5mm}21]{tanabe1992}.
    Therefore, with $p_i = \ell^{-1}$ for all $i$, we see from \eqref{np} that
    \begin{align}\label{eq:convergence.standardized.chi.squared}
        \frac{p(\bb{N}) - n \ell^{-1}}{\sqrt{n \ell^{-1} (1 - \ell^{-1})}}
        &\stackrel{\phantom{***}}{=} \sqrt{\frac{m}{m-1}} \left[\frac{\sum_{i=1}^{\ell-1} \sum_{i'=1}^{\ell-1} (\Sigma^{-1})_{i,i'} \frac{(N_i - m \ell^{-1})}{\sqrt{m}} \frac{(N_{i'} - m \ell^{-1})}{\sqrt{m}}}{\sqrt{2 (\ell - 1)}} - \sqrt{\frac{\ell - 1}{2}}  \right] \notag \\[1mm]
        &\stackrel{d}{\longrightarrow} \frac{\Gamma - (\ell - 1)}{\sqrt{2 (\ell - 1)}}, \quad \text{where } \Gamma\sim \chi_{\ell-1}^2.
    \end{align}
    Now, let
    \begin{equation}\label{definition_r_U_V_tilde}
        \widetilde{U}_k \leqdef \frac{\sigma_U}{\sigma} \cdot \frac{U_k - \mu_U}{\sigma_U} ~\quad \text{and} ~\quad \widetilde{V}_k \leqdef \frac{\sigma_V}{\sigma} \cdot \frac{V_k - \mu_V}{\sigma_V},
    \end{equation}
    then we can write
    \begin{equation}
        S_n = \frac{1}{\sqrt{n}}\Bigg(r \, \frac{\big(p(\bb{N}) - n \ell^{-1}\big)}{\sqrt{\ell^{-1} (1 - \ell^{-1})}} + \sum_{\substack{j=1 \\ D_k = 0}}^n \widetilde{U}_k + \sum_{\substack{j=1 \\ D_k = 1}}^n \widetilde{V}_k\Bigg),
    \end{equation}
    since, from \eqref{eq:F_Xk}, we know that
    \begin{equation}\label{eq:mu.decomposition}
        \mu = (1 - \ell^{-1}) \mu_U + \ell^{-1} \mu_V.
    \end{equation}
    With the notation $t_n \leqdef t / \sqrt{n}$, the mutual independence between the $U_k$'s, the $V_k$'s and $\bb{M}_{\bullet} \leqdef \{\bb{M}_j\}_{j=1}^m$ yields, for all $t\in \R$,
    \begin{align}\label{eq:thm:main.result.beginning}
        \E\big[e^{\ii t S_n} | \bb{M}_{\bullet}\big]
        &= e^{\ii t r \, \frac{(p(\bb{N}) - n \ell^{-1})}{\sqrt{n \ell^{-1} (1 - \ell^{-1})}}} \prod_{\substack{j=1 \\ D_k = 0}}^n \E[e^{\ii t_n \widetilde{U}_k} | \bb{M}_{\bullet}] \prod_{\substack{j=1 \\ D_k = 1}}^n \E[e^{\ii t_n \widetilde{V}_k} | \bb{M}_{\bullet}] \notag \\
        &= e^{\ii t r \, \frac{(p(\bb{N}) - n \ell^{-1})}{\sqrt{n \ell^{-1} (1 - \ell^{-1})}}} ~ [\varphi_{\widetilde{U}}(t_n)]^{n (1 - \ell^{-1})} [\varphi_{\widetilde{V}}(t_n)]^{n \ell^{-1}} \bigg[\frac{\varphi_{\widetilde{V}}(t_n)}{\varphi_{\widetilde{U}}(t_n)}\bigg]^{p(\bb{N}) - n \ell^{-1}} \notag \\[1mm]
        &= e^{\ii t r \, \frac{(p(\bb{N}) - n \ell^{-1})}{\sqrt{n \ell^{-1} (1 - \ell^{-1})}}} \cdot [\varphi_{\widetilde{U}}(t_n)]^{n (1 - \ell^{-1})} [\varphi_{\widetilde{V}}(t_n)]^{n \ell^{-1}} \notag \\
        &\quad\cdot \left[\frac{[\varphi_{\widetilde{V}}(t_n)]^n \cdot e^{\frac{1}{2} \cdot \frac{\sigma_V^2}{\sigma^2} t^2}}{[\varphi_{\widetilde{U}}(t_n)]^n \cdot e^{\frac{1}{2} \cdot \frac{\sigma_U^2}{\sigma^2} t^2}}\right]^{\frac{p(\bb{N}) - n \ell^{-1}}{n}} \hspace{-6mm} \cdot \hspace{4mm} \left[\frac{e^{-\frac{1}{2} \cdot \frac{\sigma_V^2}{\sigma^2} t^2}}{e^{-\frac{1}{2} \cdot \frac{\sigma_U^2}{\sigma^2} t^2}}\right]^{\frac{p(\bb{N}) - n \ell^{-1}}{n}}\hspace{-6mm}.
    \end{align}
    (The reader should note that, for $n$ large enough, the manipulations of exponents in the second and third equality above are valid because the highest powers of the complex numbers involved have their principal argument converging to $0$. This stems from the fact that $0 \leq p(\bb{N}) \leq n$, and the quantities $[\varphi_{\widetilde{U}}(t_n)]^n$ and $[\varphi_{\widetilde{V}}(t_n)]^n$ both converge to real exponentials as $n\to \infty$, by the CLT.)
    We now evaluate the four terms on the right-hand side of \eqref{eq:thm:main.result.beginning}.
    For the first term in \eqref{eq:thm:main.result.beginning}, the continuous mapping theorem and \eqref{eq:convergence.standardized.chi.squared} yield
    \begin{equation}\label{eq:thm:main.result.beginning.part.1}
        e^{\ii t r \, \frac{(p(\bb{N}) - n \ell^{-1})}{\sqrt{n \ell^{-1} (1 - \ell^{-1})}}} \stackrel{d}{\longrightarrow} e^{\ii t r \frac{\Gamma - (\ell - 1)}{\sqrt{2 (\ell - 1)}}}, \quad \text{as } m\to \infty.
    \end{equation}
    For the second term in \eqref{eq:thm:main.result.beginning}, the CLT yields
    \begin{equation}\label{eq:thm:main.result.beginning.part.2}
        \begin{aligned}
            [\varphi_{\widetilde{U}}(t_n)]^{n (1 - \ell^{-1})} [\varphi_{\widetilde{V}}(t_n)]^{n \ell^{-1}}
            &\underset{m\to \infty}{\longrightarrow} \exp\Big(-\frac{1}{2} \cdot (1 - \ell^{-1}) \frac{\sigma_U^2}{\sigma^2} t^2\Big) \exp\Big(-\frac{1}{2} \cdot \ell^{-1} \frac{\sigma_V^2}{\sigma^2} t^2\Big) \\[2mm]
            &=e^{-\frac{1}{2} (1 - r^2) t^2},
        \end{aligned}
    \end{equation}
    where in the last equality we used that, from \eqref{eq:F_Xk}, we have
    \begin{align}
        \sigma^2
        &= \E[X^2] - \mu^2 = (1 - \ell^{-1}) \sigma_{U}^2 + \ell^{-1} \sigma_{V}^2 + \ell^{-1} (1 - \ell^{-1}) (\mu_U - \mu_V)^2. \label{eq:relation.r:var}
    \end{align}
    For the third term in \eqref{eq:thm:main.result.beginning}, the quantity inside the bracket converges to $1$ by the CLT.
    Hence, the elementary bound
    \begin{equation}
        |e^z - 1| \leq |z| + \sum_{j=2}^{\infty} \frac{|z|^j}{2} \leq |z| + \frac{|z|^2}{2 (1 - |z|)} \leq \frac{1 + \ell^{-1}}{2 \ell^{-1}} |z|, ~\quad \text{for all } |z| \leq 1 - \ell^{-1},
    \end{equation}
    and the fact that $\big|\frac{p(\bb{N}) - n \ell^{-1}}{n}\big| \leq 1 - \ell^{-1}$ yield, as $m\to \infty$,
    \begin{equation}\label{eq:thm:main.result.beginning.part.3}
        \left|\left[\frac{[\varphi_{\widetilde{V}}(t_n)]^n \cdot e^{\frac{1}{2} \cdot \frac{\sigma_V^2}{\sigma^2} t^2}}{[\varphi_{\widetilde{U}}(t_n)]^n \cdot e^{\frac{1}{2} \cdot \frac{\sigma_U^2}{\sigma^2} t^2}}\right]^{\frac{p(\bb{N}) - n \ell^{-1}}{n}} - 1\right| \leq \frac{1 - \ell^{-2}}{2 \ell^{-1}} \left|\text{Log}\left[\frac{[\varphi_{\widetilde{V}}(t_n)]^n \cdot e^{\frac{1}{2} \cdot \frac{\sigma_V^2}{\sigma^2} t^2}}{[\varphi_{\widetilde{U}}(t_n)]^n \cdot e^{\frac{1}{2} \cdot \frac{\sigma_U^2}{\sigma^2} t^2}}\right]\right| \longrightarrow 0.
    \end{equation}
    For the fourth term in \eqref{eq:thm:main.result.beginning}, the continuous mapping theorem and $\frac{p(\bb{N}) - n \ell^{-1}}{n} \stackrel{\PP}{\longrightarrow} 0$ (recall \eqref{eq:convergence.standardized.chi.squared}) yield
    \begin{equation}\label{eq:thm:main.result.beginning.part.4}
        \left[\frac{e^{-\frac{1}{2} \cdot \frac{\sigma_V^2}{\sigma^2} t^2}}{e^{-\frac{1}{2} \cdot \frac{\sigma_U^2}{\sigma^2} t^2}}\right]^{\frac{p(\bb{N}) - n \ell^{-1}}{n}} \stackrel{\PP}{\longrightarrow} 1, \quad \text{as } m\to \infty.
    \end{equation}
    By combining \eqref{eq:thm:main.result.beginning.part.1}, \eqref{eq:thm:main.result.beginning.part.2}, \eqref{eq:thm:main.result.beginning.part.3} and \eqref{eq:thm:main.result.beginning.part.4}, Slutsky's theorem implies, for all $t\in \R$,
    \begin{equation}\label{eq:thm:main.result.beginning.combined}
        \E\big[e^{\ii t S_n} | \bb{M}_{\bullet}\big] \stackrel{d}{\longrightarrow} e^{\ii t r \frac{\Gamma - (\ell - 1)}{\sqrt{2 (\ell - 1)}}} e^{-\frac{1}{2} (1 - r^2) t^2}, \quad \text{as } m\to \infty.
    \end{equation}
    Since the sequence $\{|\E[e^{\ii t S_n} | \bb{M}_{\cdot}]|\}_{m\in \N}$ is uniformly integrable (it is bounded by $1$), Theorem 25.12 in \cite{billingsley1995} shows that we also have the mean convergence
    \begin{equation}\label{eq:thm:main.result.mean.convergence}
        \E\big[\E\big[e^{\ii t S_n} | \bb{M}_{\bullet}\big]\big] \longrightarrow \E\Big[e^{\ii t r \frac{\Gamma - (\ell - 1)}{\sqrt{2 (\ell - 1)}}}\Big] e^{-\frac{1}{2} (1 - r^2) t^2}, \quad \text{as } m\to \infty,
    \end{equation}
    which proves \eqref{cf_S}. The conclusion follows.\qed
\end{proofff}

\section{Properties of \textit{S}} \label{properties-of-s}

    Recall that $F$ denotes the marginal distribution of the r.v.s $X_1, \ldots, X_n$ in~\eqref{the_X_sequence}. Theorem \ref{limiting_Sn_thm} states that $S_n$, the standardized sample mean of these r.v.s, converges in distribution to a r.v.\ $S$ whose characteristic function is given by~\eqref{cf_S}.
    When the choice $\ell \geq 2$ is fixed, the distribution of $S$ has only one parameter, $r$, defined as
    \begin{equation}\label{r_definition}
        r = \frac{\sqrt{\ell^{-1} (1 - \ell^{-1})} (\mu_V - \mu_U)}{\sigma} = \frac{\mu_V - \mu}{\sigma \sqrt{\ell - 1}},
    \end{equation}
    where the second equality stems from \eqref{eq:mu.decomposition}. Hence, $r$ depends on the margin $F$ (through the quantities $A$, $\mu_U$, $\mu_V$ and $\sigma$).
    The behavior of $S$ with respect to $F$ (via $r$) is now studied. From \eqref{eq:relation.r:var}, we see that
    \begin{equation}\label{bounds_on_r}
        r^2 = 1 - \frac{(1 - \ell^{-1}) \sigma_U^2 + \ell^{-1} \sigma_V^2}{\sigma^2}, \quad \text{and thus} \quad 0 \leq r^2 \leq 1.
    \end{equation}
        Example \ref{clt_logNorm_example} and several other examples in Appendix~\ref{appendix_examples} show how the critical points $r^2 = 0, 1$ can be achieved or approached when $F$ is discrete or absolutely continuous. See also  Appendix~\ref{appendix_computing_code} for the \texttt{R} computing codes to generate observations from all these examples. (Note that these examples could serve as scenarios of dependence to compare, via Monte-Carlo simulations, various tests of independence; see e.g., \cite{huskova2008}.)

     \begin{examp}[$r$ arbitrarily close to $0$ when $F$ is absolutely continuous]\label{clt_logNorm_example}
        Let $\ell = 2$, $A = [1,\infty)$ and let $W \sim f $ where $f$ is the density of a Log-Normal($0, \beta)$ distribution. Note that $\textup{median}(W) = 1$, $\E[W] = \exp(\beta/2)$, and $\Var[W] = \left[\exp(\beta) - 1 \right] \exp(\beta)$. Furthermore,
        \begin{equation}
            \mu_V = \int_{1}^{\infty} 2 x f(x) dx = \sqrt{\frac{2}{\pi \beta}} \int_{1}^{\infty}  \exp\left(-\frac{(\log x)^2}{2 \beta} \right) \rd x =  \exp(\beta/2) \left[1 + \textup{Erf}\big(\sqrt{\beta / 2}\big)\right],
        \end{equation}
        where the integral on the second line was solved with \texttt{Mathematica}. Hence, from~\eqref{r_definition},
        \begin{equation}
            r = \frac{\mu_V - \E[W]}{\sqrt{\Var[W]}}  = \frac{\exp(\beta/2) \textup{Erf}\big(\sqrt{\beta / 2}\big)}{\sqrt{\left[\exp(\beta) - 1 \right] \exp(\beta)}} =  \frac{\textup{Erf}\big(\sqrt{\beta / 2}\big)}{\sqrt{\exp(\beta) -1}},
        \end{equation}
        and it is straightforward to see that $r \to 0$ as $\beta \to \infty$.
    \end{examp}

\begin{figure}[ht]
        \centering
        \hspace{-5mm}
        \includegraphics[scale=0.35]{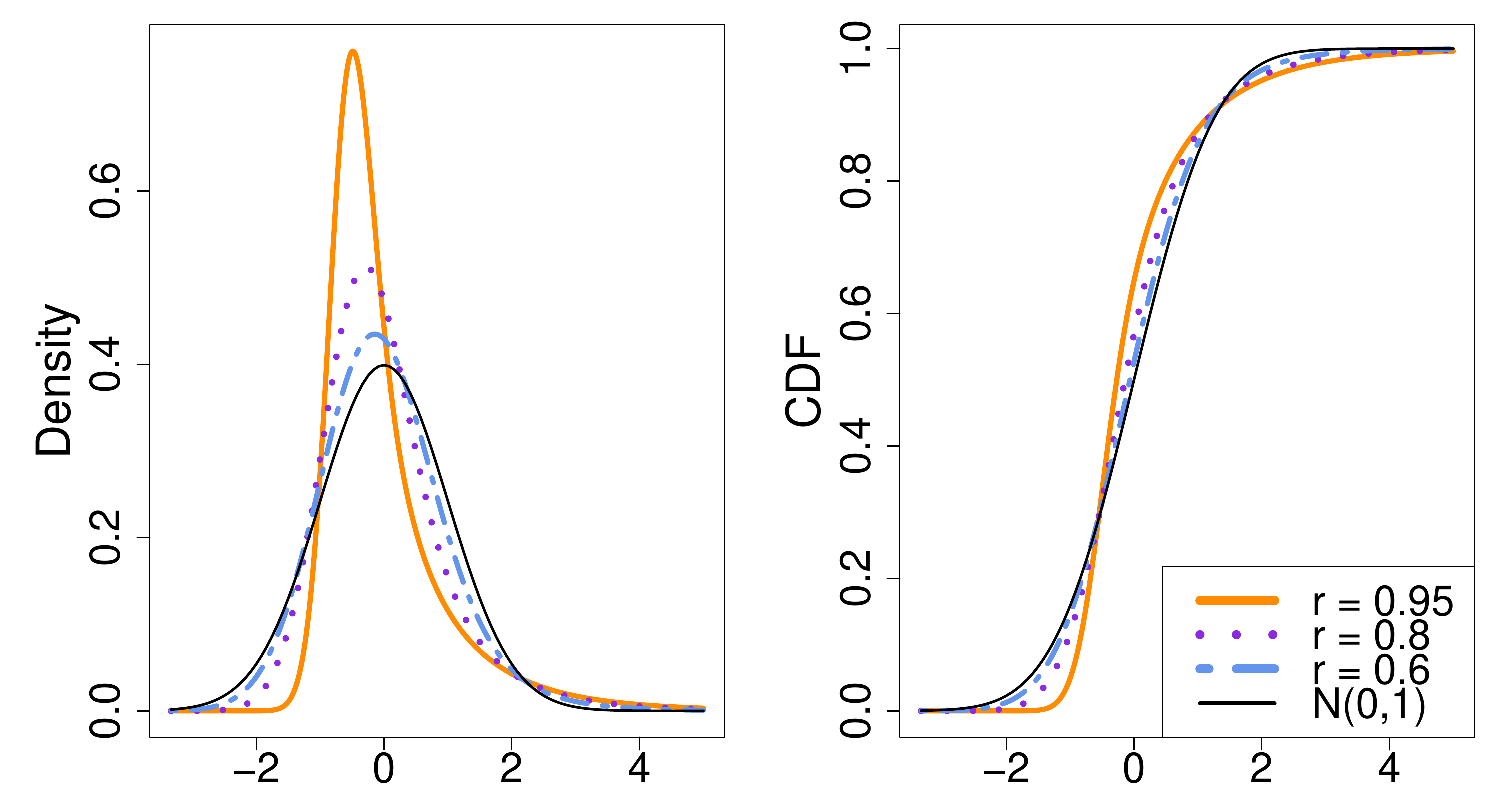}
        \captionsetup{width=0.85\textwidth}
        \caption{Density (left) and CDF (right)  of $S$ for fixed $\ell = 2$ and varying $r$ ($r=0.6, 0.8, 0.95)$, compared to those of a $N(0,1)$. This illustrates that the CLT can `fail' substantially under pairwise independence.}
        \label{cdf_df_r2_comparison}
    \end{figure}

\begin{figure}[ht]
        \centering
        \hspace{-5mm}
        \includegraphics[scale=0.35]{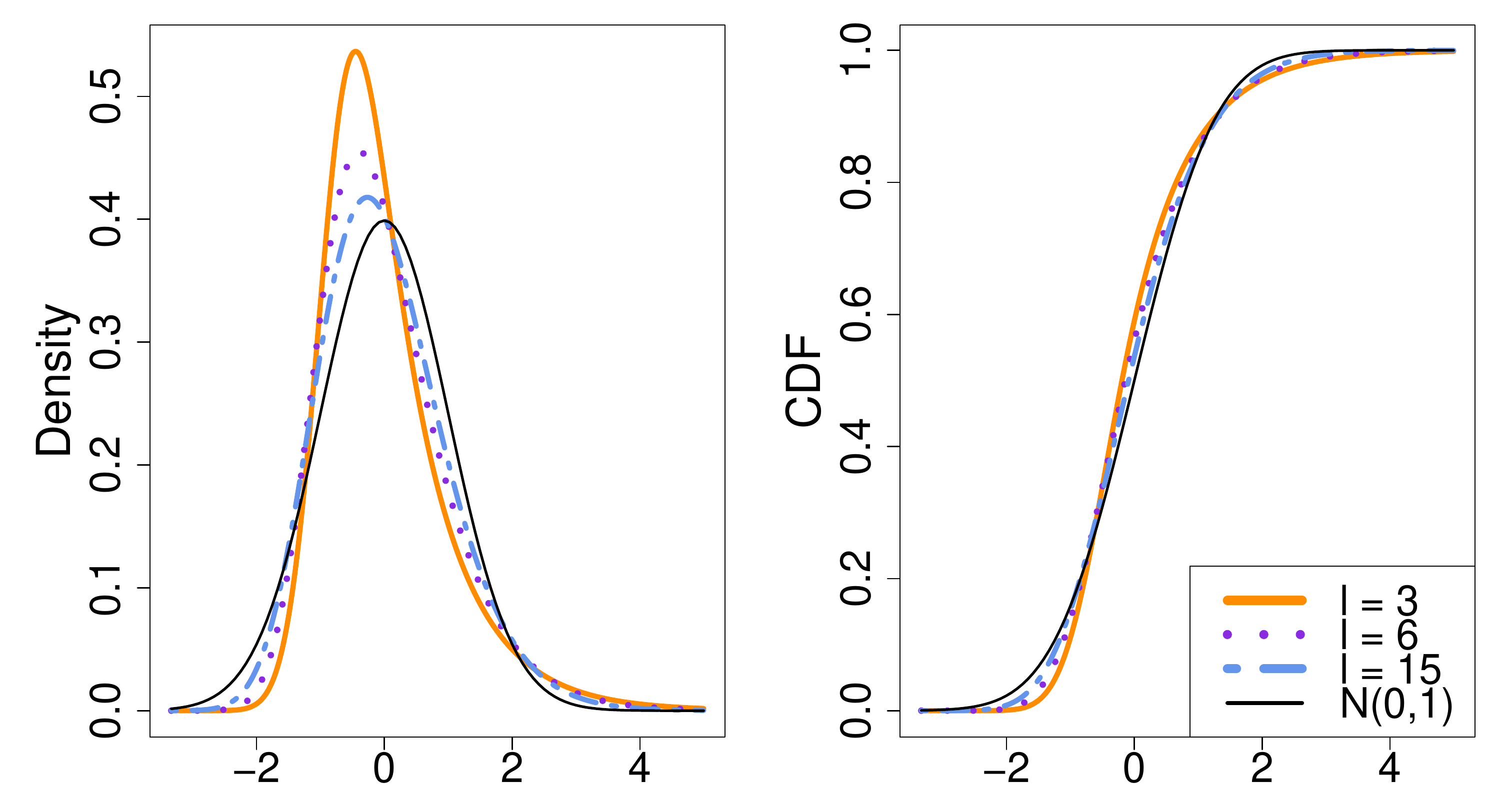}
        \captionsetup{width=0.85\textwidth}
        \caption{Density (left) and CDF (right) of $S$ for fixed $r = 0.9$ and varying $\ell$ ($\ell=3, 6, 15$), compared to those of a $N(0,1)$. This illustrates that $S$ converges to a $N(0,1)$ as $\ell$ grows.}
        \label{cdf_df_ell_comparison}
    \end{figure}

    Next, recall that the characteristic function on the right-hand side of \eqref{cf_S} is that of
    \begin{equation}
        S = \sqrt{1 - r^2} Z + r \, \chi,
    \end{equation}
    where the r.v.s $Z\sim N(0,1)$ and $\chi\sim [\chi_{\ell-1}^2 - (\ell - 1)]/\sqrt{2 (\ell - 1)}$ are independent. This makes it clear that, when $\ell \geq 2$ is fixed, $r$ completely determines the shape of $S$; the closer $r$ gets to $0$, the closer the distribution of $S$ is to a standard Gaussian, while the closer $r$ gets to $\pm 1$, the closer the distribution of $S$ is to a standardized $\pm \chi^2_{\ell-1}$.
    This shift from a Gaussian distribution towards a $\chi^2_{\ell-1}$ distribution is represented graphically in Figure~\ref{cdf_df_r2_comparison} (where $\ell = 2$ and $r$ varies). On the other hand, regardless of $r$, if $\ell$ increases then $S$ gets closer to a $N(0,1)$, as illustrated in Figure \ref{cdf_df_ell_comparison} (where $r=0.9$ and $\ell$ varies). These figures illustrate clearly that pairwise independence might be a very poor substitute to mutual independence as an assumption in the CLT.

    In terms of moments, simple calculations with \texttt{Mathematica} yield that
    \begin{equation}
        \E[S] = 0, \quad \E[S^2] = 1, \quad \E[S^3] = \sqrt{\frac{8}{\ell - 1}} \, r^3 \quad \text{and} \quad \E[S^4] = 3 + \frac{12}{\ell - 1} r^4,
    \end{equation}
    so that upper bounds on the skewness and kurtosis of $S$ are $\sqrt{8 / (\ell - 1)}$ and $3 + 12 / (\ell - 1)$, respectively. The limiting r.v.\ $S$ can therefore be much more skewed and heavy-tailed than the standard Gaussian distribution, which is also confirmed by Figure~\ref{cdf_df_r2_comparison}.  


    Lastly, let us comment on the $r$ parameter, and explain why a $r$ close to $1$ yields a more `drastic' failure of the CLT. First, recall that the CLT fails when applied to the sequence of pairwise independent r.v.s $D_1, \ldots, D_n$ given in~\eqref{the_D_sequence} because the proportion of $1$'s in that sequence can be very large, whereas the proportion of $0$'s can never be large. Consequently, the distribution of the asymptotic sample mean of this sequence is asymmetrical (skewed to the right). When we `assign' an arbitrary margin to the $D$'s in order to create our sequence $\{X_j , j \geq 1\}$, we can attenuate (to a certain degree) this asymmetry.
    Consider for example the case $\ell = 2$ and $A = [\widetilde{w},\infty)$, where $\widetilde{w}$ denotes the median of an absolutely continuous distribution $F$. In this case, the $X$'s, as opposed to the $D$'s, take a continuous range of values, and hence $X$'s `above the median' can be quite close to their mean (whereas the $D$'s are all either `much bigger' or `much smaller' than their mean). The parameter $r = (\mu_V - \mu)/\sigma$ measures to what extent this `attenuation of asymmetry' happens. 
    Indeed, if $r$ is close to $0$, the $X$'s observations above the median are not too far away from the mean (on average). This implies that, even if the proportion of observations above the median is huge, it will not overly boost the overall mean of the sample, and the distribution of this mean will not be overly asymmetrical.

    To give a concrete example (again with $\ell =2, A = [\widetilde{w}, \infty))$, let
    $X\sim$ Log-normal($\alpha, \beta$). In that case, simple calculations (see Example~\ref{clt_logNorm_example} for details) yield $r = \textup{Erf}\big(\sqrt{\beta/2}\big)/ \sqrt{\exp(\beta) -1}$, a decreasing function of~$\beta$. On the other hand, it is well known that the kurtosis of $X$ is an increasing function of $\beta$. So, increasing $\beta$ makes $X$ heavier tailed, while giving a lower value of $r$. For ease of interpretation, and since $r$ is invariant to shifting and scaling,  consider the r.v.\ $Y = (X-\mu)/\sigma$, which has the same value of $r$ and the same kurtosis as $X$. For $Y$, it is clear from \eqref{r_definition} that $r$ is just the \textit{mean of $Y$ given that it exceeds its median} (i.e., $\E[Y|Y>\textrm{median}(Y)]$). As $\beta$ increases, the right tail gets longer. To compensate the more extreme values on the right, while keeping the mean of $Y$ equal to $0$, the median is forced to move further away to the left of the mean. Hence, the mean of observations above the median, i.e.\ $r$, also gets smaller. 





\section{Conclusion}\label{conclusion_section}

    We showed that the CLT can `fail' for a pairwise independent sequence of identically distributed r.v.s $\{X_j, j \geq 1\}$ having \textit{any} distribution that satisfies Condition \ref{cond:F}. Under a specific structure of dependence for such pairwise independent $X_j$'s, we obtained the asymptotic distribution of the standardized sample mean $S_n$ and found it to be always `worse behaved' than a Gaussian. Furthermore, the extent of this departure from normality depends on the initial common margin of the $X_j$'s. 
    This is in contradiction with the CLT under which, regardless of the margin, $S_n$ always converges to a Gaussian.

    A corollary of our main result is that there exists a sequence of pairwise independent Gaussian r.v.s for which the limiting distribution of $S_n$ is substantially `worse behaved' than a Gaussian, being asymmetric and heavier tailed. To our knowledge, no other such example exists in the literature. Given the widespread use of the CLT, even in standard parametric statistical techniques such as tests and confidence intervals for means and variances  \citep[see, e.g.,][]{Coeurjolly2009a}, this constitutes 
    a serious warning to practitioners of statistics who may think that, to invoke the CLT, all one needs is for the original random variables $\{X_j, j \geq 1\}$ to be approximately Gaussian or to have a large enough sample size. Mutual independence \textit{is} a crucial assumption that should not be forgotten, nor misunderstood. 



    As a final note, our sequence is not strictly stationary, and it is not obvious that there exists a stationary sequence with a similar asymptotic distribution for $S_n$. Furthermore, some authors have studied the CLT under $K$-tuplewise independence (for $K > 2$); see, e.g., \citet{pruss1998, bradley2009, bradley2010, weakley2013}. It would be interesting to generalize our construction in that direction. One might wonder if, as $K$ increases, the distribution of $S_n$ would necessarily get closer to that of a Gaussian. An articulated answer to this question is not trivial and is left for future research.

\appendix

\begin{appendices}

\section{Other examples}\label{appendix_examples}

    \begin{examp}[$r = 0$ when $F$ is discrete]
        For any integer $\ell \geq 2$, take
        \begin{equation}
            A = \{-1,1\}, \quad \PP(W = -1) = \PP(W = 1) = \tfrac{\ell^{-1}}{2} \quad \text{and} \quad \PP(W = -2) = \PP(W = 2) = \tfrac{1 - \ell^{-1}}{2},
        \end{equation}
        since it implies $\mu_U = \mu_V = 0$.
    \end{examp}

    \begin{examp}[$r = 0$ when $F$ is absolutely continuous]
        For any integer $\ell \geq 2$, take
        \begin{equation}
            A = [-\ell^{-1}, \ell^{-1}] \quad \text{and} \quad W\sim \text{Uniform}\hspace{0.3mm}[-1,1],
        \end{equation}
        since again it implies $\mu_U = \mu_V = 0$.
    \end{examp}

    \begin{examp}[$r^2 = 1$ when $F$ is discrete]
        Let $\ell \geq 2$ be any integer.
        To get $r = 1$, take
        \begin{equation}
            A = \{1\}, \quad \PP(W = 1) = \ell^{-1} \quad \text{and} \quad \PP(W = -1) = 1 - \ell^{-1},
        \end{equation}
        since this means $\sigma_U = \sigma_V = 0$ and $\mu_V > \mu_U$. By symmetry, taking $A^c = \{1\}$ instead yields $r = -1$.
    \end{examp}

    \begin{examp}[$r$ arbitrarily close to $1$ when $F$ is absolutely continuous]\label{ex_r_close_1}
        Let $f_1$ be the density function of a N($-\ell^{-1}, \sigma^2$), $f_2$ the density function of a N($1 - \ell^{-1}, \sigma^2$), and $f$ their mixture: $f(x) = (1 - \ell^{-1}) f_1(x) + \ell^{-1} f_2(x)$. Then, for $W \sim f$, we have $\E[W] = 0$ and $\Var[W] = \E[W^2] = \sigma^2 + \ell^{-1} (1 - \ell^{-1})$.
        Assuming that $0 < \sigma \leq \tfrac{1}{2} \ell^{-1}$, a straightforward Gaussian tail estimate on $f_1$ shows that there exists $w_{\ell} \in (-\ell^{-1},1 - \ell^{-1})$ such that $\PP(W\in [w_{\ell},\infty)) = \ell^{-1}$.
        If we take $A = [w_{\ell},\infty)$, then we have
        \begin{align}
            \mu_V
            &= \ell \int_{w_{\ell}}^{\infty} x \, ((1 - \ell^{-1}) f_1(x) + \ell^{-1} f_2(x)) \, \rd x \notag \\
            &= (\ell - 1) \int_{w_{\ell}}^{\infty} \Big(\frac{x + \ell^{-1}}{\sigma} - \frac{\ell^{-1}}{\sigma}\Big) \frac{1}{\sqrt{2\pi}} \exp\left(-\frac{1}{2} \Big(\frac{x + \ell^{-1}}{\sigma}\Big)^2\right) \rd x \notag \\
            &\quad+ \int_{w_{\ell}}^{\infty} \Big(\frac{x - (1 - \ell^{-1})}{\sigma} + \frac{1 - \ell^{-1}}{\sigma}\Big) \frac{1}{\sqrt{2\pi}} \exp\left(-\frac{1}{2} \Big(\frac{x - (1 - \ell^{-1})}{\sigma}\Big)^2\right) \rd x \notag \\
            &= (\ell - 1) \frac{\sigma}{\sqrt{2\pi}} \exp\left(-\frac{1}{2} \Big(\frac{w_{\ell} + \ell^{-1}}{\sigma}\Big)^2\right) - (1 - \ell^{-1}) \Psi\Big(\frac{w_{\ell} + \ell^{-1}}{\sigma}\Big) \notag \\
            &\quad+ \frac{\sigma}{\sqrt{2\pi}} \exp\left(-\frac{1}{2} \Big(\frac{w_{\ell} - (1 - \ell^{-1})}{\sigma}\Big)^2\right) + (1 - \ell^{-1}) \Psi\Big(\frac{w_{\ell} - (1 - \ell^{-1})}{\sigma}\Big),
        \end{align}
        where $\Psi(z)$ is the survival function of the standard Gaussian.
        Therefore, from \eqref{r_definition}, we have $r\to 1$ as $\sigma\to 0$ :
        \begin{equation}
            r = \frac{\mu_V - \E[W]}{\sqrt{\Var(W)} \sqrt{\ell - 1}} \underset{\sigma\to 0}{\longrightarrow} \frac{0 - 0 + 0 + (1 - \ell^{-1}) - 0}{\sqrt{0^2 + \ell^{-1} (1 - \ell^{-1})} \sqrt{\ell - 1}} = 1.
        \end{equation}
    \end{examp}

    \begin{examp}[$F$ is a $N(\mu, \sigma^2)$]
        Let $\ell = 2$, choose $A = [\mu,\infty)$ and let $Z$ be a $N(0,1)$ r.v. Then,
        \begin{equation}
           \mu_V = \mu + \sigma \, \E[Z|Z>0] = \mu + \sigma \, \E[|Z|] = \mu + \sigma \, \sqrt{\frac{2}{\pi}}.
        \end{equation}
        It follows that $r = \sqrt{2/\pi} \approx 0.8$ (irrespective of $\mu$ and $\sigma$). Note that this corresponds to the purple dotted curve on Figure \ref{cdf_df_r2_comparison}. Hence this case provides a nice illustration of how `badly' the CLT can fail for pairwise independent Gaussian variables.
    \end{examp}

\section{Computing codes}\label{appendix_computing_code}

\begin{lstlisting}[language=R]
# Generator of pairwise independent observations
piid.generator <- function(m = 3, randF = rnorm, 
                           indA = function(x) ifelse(x <= 0, FALSE, TRUE),
                           ell = 2) {
  
  # Check that the value of 'ell' provided is coherent 
  #  with the function 'indA' provided
  # This check is valid only for moderate values of 'ell'
  if (round(1 / mean(indA(randF(10 ^ 5)))) != ell) 
    warning("Is 'ell' consistent with your A?")
  
  # Sample size
  n <- choose(m, 2)
  
  # Generate the 'initial' multinomial sample of size m
  M <- rmultinom(m, 1, rep(1 / ell, ell))
  
  # Find all possible pairs out of the m multinomials
  # Use those pairs to create the 'n' D's in Eq. (2.4)
  combin <- combn(1:m, 2)
  D <- apply(M[, combin[1,]] == M[, combin[2,]], 2, all)
  D <- as.integer(D)
  # Compute the number of 1's among the D's
  pN <- sum(D)
  
  # Generate pN r.v.s with distribution F restricted to A 
  # and (n - pN) r.v.s with distribution F restricted to A^c
  nU <- 0
  nV <- 0
  U <- rep(NA, n - pN)
  V <- rep(NA, pN)
  while((nU < n - pN) | (nV < pN)) {
    W <- randF(1)
    indAW <- indA(W)
    if (indAW & (nV < pN)) {
      nV <- nV + 1
      V[nV] <- W
    } else if ((indAW == 0) & (nU < n - pN)) {
      nU <- nU + 1
      U[nU] <- W
    }
  }
  # Return the resulting random generated variables  
  X <- rep(NA, n)
  X[which(D == 0)] <- U
  X[which(D == 1)] <- V
  return(X)
}

\end{lstlisting}

\begin{lstlisting}[language=R]
# Specific values of randF, indA and ell for all examples

# Example 5
randF <- function(m, beta) rlnorm(m, 0, beta)
indA  <- function(x) ifelse(x >= 1, TRUE, FALSE)
piid.generator(randF = function(m) randF(m, beta = 2), indA = indA, ell = 2L)

# Example 6
r.ex6 <- function(rand, ell) {
  if (rand < 1 / (2 * ell)) {
    res <- -1L
  } else if (rand < 1 / ell) {
    res <- 1L
  } else if (rand < 1 / (2 * ell) + 1 / 2) {
    res <- -2L
  } else {
    res <- 2L
  }
  return(res)
}
randF <- function(m, ell) sapply(runif(m), FUN = r.ex6, ell = ell)
indA  <- function(x) ifelse((x == 1L) | (x == -1L), TRUE, FALSE)
piid.generator(randF = function(m) randF(m, ell = 2L), indA = indA, ell = 2L)

# Example 7
randF <- function(m) runif(m, -1, 1)
indA  <- function(x, ell) ifelse((x >= -1 / ell) & (x <= 1 / ell), TRUE, FALSE)
piid.generator(randF = randF, indA = function(x) indA(x, ell = 2L), ell = 2L)

# Example 8
randF <- function(m, ell) as.integer(2 * rbinom(m, 1, 1 / ell) - 1)
indA  <- function(x) ifelse(x == 1L, TRUE, FALSE)
piid.generator(randF = function(m) randF(m, ell = 10L), indA = indA, ell = 10L)

# Example 9
F  <-  function(x, ell, sigma) sum(c(1 - 1 / ell, 1 / ell) *
                         pnorm(x, mean = c(-1 / ell, 1 - 1 / ell),
                               sd = c(sigma, sigma)))
Finv  <- function(p, ell, sigma = 1/(4*ell)){ 
  G = function(x) F(x, ell, sigma) - p
  return(uniroot(G, c(-100,100))$root)
}
r.ex9 <- function(rand, ell, sigma) {
  if (rand < 1 / ell) {
    res <- rnorm(1, 1 - 1 / ell, sigma)
  } else {
    res <- rnorm(1, -1 / ell, sigma)
  }
  return(res)
}
randF <- function(m, ell, sigma = 1/(4*ell)) 
          sapply(runif(m), FUN = r.ex9, ell = ell, sigma = sigma)
indA  <- function(x, ell, sigma = 1/(4*ell)) 
          ifelse(x >= Finv(1 - 1 / ell, ell = ell, sigma = sigma), TRUE, FALSE)
piid.generator(randF = function(m) randF(m, ell = 3L), 
                indA = function(x) indA(x, ell = 3L), ell = 3L)

# Example 10
randF <- function(m, mu, sigma) rnorm(m, mu, sigma)
indA  <- function(x, mu) ifelse(x >= mu, TRUE, FALSE)
piid.generator(randF = function(m) randF(m, mu = 2, sigma = 1),
               indA = function(x) indA(x, mu = 2), ell = 2L)
\end{lstlisting}

\end{appendices}

\section*{Acknowledgements}

\noindent
Elements of this paper were presented at the conference Perspectives on Actuarial Risks in Talks of Young Researchers (Sibiu, Romania) in April 2019, at the $23^\text{nd}$ International Congress on Insurance: Mathematics and Economics (Munich, Germany) in July 2019 and at the Australasian Actuarial Education and Research Symposium (Melbourne, Australia) in November 2019. The authors are grateful for constructive comments received from colleagues at those conferences.

\section*{Conflict of interest}
\noindent
The authors have no conflict of interest to disclose.


%
%

\bibliographystyle{authordate1}
\bibliography{pairwise_indep_non_CLT_bib}

\end{document}